\documentclass[12pt]{article}
\title{Fibonomial Cumulative Connection Constants}
\author{A.K.Kwa\'sniewski \\
Member of the Institute of Combinatorics and its Applications\\ 
Higher School of Mathematic and Applied Informatics\\
PL-15-021 Bia{\l}ystok, ul.Kamienna 17, POLAND\\
e-mail: kwandr@gmail.com}

\usepackage{amsmath,amsthm}
\usepackage{graphicx}
\usepackage{color}
\chardef\bslash=`\\ % p. 424, TeXbook
\newtheorem{defn}{Definition}[section]
\newtheorem{obs}{Observation}[section]

\newtheorem{rem}{Remark}[section]

\newcommand{\fnomial}[2]{ {{#1} \choose {#2}}_F }

\begin{document}
\maketitle
%%%%%%%%%%%%
\begin{abstract}
\noindent In this note we present examples of cumulative connection
constants  - new Fibonomial  ones included. All examples posses
combinatorial interpretation.\\

\noindent This presentation is an update  of A.K.Kwa\'sniewski  \textit{Fibonomial cumulative connection constants}  \textbf{Bulletin of the ICA vol. 44 (2005), 81-92.}
\end{abstract}

\vspace{0.2cm}

\noindent AMS Classification Numbers: 06A06 ,05B20, 05C75, 05A40

\vspace{0.2cm}

\noindent this electronic paper is affiliated to The Internet Gian-Carlo Polish Seminar:

\noindent \emph{http://ii.uwb.edu.pl/akk/sem/sem\_rota.htm}

\vspace{0.2cm}

\section{Introduction}

\noindent The cumulative connection constants  ({\bf ccc}) $C_{n}$ were introduced
  by the Authors of  \cite{1}:
  $$ C_{n}=\sum_{k\geq 0}c_{n,k}\; ,\;\;\;n\geq 0,$$
  where
$$p_{n}(x)=\sum_{0\leq k \leq n}c_{n,k}q_{k}(x) $$
   and Pascal-like array $\left\{ c_{n,k}\right\}_{n,k\geq 0}$  of connection constants
    $c_{n,k},\;\;n,k\geq 0$ "connects" two polynomial sequences  (Note: $\left\{ p_{k}(x)\right\}_{k\geq 0}$ is a polynomial
sequence if $\deg p_{k} = k $).

\noindent \underline{Motivation 1}\\
\vspace{0.2cm}
\noindent  "The connection constants problem"
i.e. combinatorial interpretations, algorithms of calculation ,
recurrences  and other properties - is  one of  the
\underline{central} issues of the binomial enumeration in Finite
Operator Calculus (FOC ) of Rota-Roman and Others and in  its
afterwards extensions (see: abundant references in
\cite{8}-\cite{10},\cite{13}). Important knowledge  on recurrences
for $\left\{ c_{n,k}\right\}_{n,k\geq 0}$ one draws from \cite{2}
- consult also the NAVIMA group program:\\
http://webs.uvigo.es/t10/navima/ .\\
\\
\noindent \underline{Motivation 2}\\
\vspace{0.1cm}
\noindent The cumulative connection constants
({\bf ccc}) as defined above appear to represent combinatorial
quantities of primary importance  - as shown by examples below .

\vspace{0.1cm}
\noindent A special case of interest  considered in \cite{2,1}  is the case
of monic persistent root polynomial  sequences \cite{3}
characterized by the following conditions:
$$p_{0}(x)=q_{0}(x)=1\; ;\;\;q_{n}(x)=q_{n-1}(x)(x-r_{n})\; ,\;p_{n}(x)=p_{n-1}(x)(x-s_{n})\;\;\;\;n\geq 1.$$

\vspace{0.1cm}
\noindent The pair of persistent root polynomial sequences and hence
corresponding connection constants array  are then bijectively
labeled by the pair of root sequences :
$$[r]=\left\{ r_{k}\right\}_{k\geq 1},\;\;\;[s]=\left\{ s_{k}\right\}_{k\geq 1}.$$
In this particular case of persistent root polynomial sequences
the authors of  \cite{2} gave a particular name to connection
constants : "generalized Lah numbers" -  now denoted as follows
$c_{n,k}=L_{n,k}$. The generalized Lah numbers do satisfy  the
following recurrence \cite{2,1} :
\begin{equation}\label{L}
L_{n+1,k}=L_{n,k-1}+(r_{k+1}-s_{n+1})L_{n,k}
\end{equation}
$$L_{0,0}=1\; , \;\;L_{0,-1}=0\;\;\;n,k\geq 0$$
from which  the recurrence for cumulative connection constants
({\bf ccc}) follows
\begin{equation}\label{C}
C_{n+1}=(1-s_{n+1})C_{n}+\sum_{k=0}^{n}c_{n,k}r_{k+1}
\end{equation}
where $c_{n,k}=L_{n,k}$.\\

 The clue examples of  \cite{1}  are given by the following choice of a pair of  root sequences
 :\\
\underline{the Fibonacci sequence}\\
\textrm{}\\ $ C_{n} = F_{n} ;\;\;  n>0\; \longleftrightarrow
[r]=\left\{ r_{k}\right\}_{k\geq
1}=\left\{0,0,0,0,0,0,\ldots\right\},\; [s]=\left\{
s_{k}\right\}_{k\geq 1}=\left\{0,0,...,0,...\right\}=
[0]$\\
\\
\underline{the Lucas sequence}\\
\textrm{}\\ $ C_{n} = L_{n}\; ;\; n>0 \longleftrightarrow
[r]=\left\{ r_{k}\right\}_{k\geq 1};\;\;   r_{1} = 0 , r_{2} = 2 ,
r_{k} = 1- (r_{k-1})^{2}$ for $k>2$ and $[s] = [0]$.

\section{Examples of cumulative connection constants - old and new} 

\vspace{0.1cm}
\noindent We shall supply now several examples of connection constants arrays and corresponding
ccc`s  - Fibonomial case included - altogether with  their
combinatorial interpretation. (HELP! - For notation-help  - see:
Appendix ).\\

\noindent {\bf Example 1.} $ c_{n,k}=\binom{n}{k},\;\;x^{n}=\sum_{k\geq
0}\binom{n}{n}(x-1)^{k}\;\;\;n\geq
0;\\E^{-1}(D)x^{k}=(x-1)^{k}\;\;n\geq 0.$
$$C_{n}=\sum_{k\geq 0}\binom{n}{k}=2^{n}= number\;\; of\; all\;  subsets\; of\;   S;\;\;   n = |S|;   n\geq 0 .$$

\noindent Here  $E^{a}(D)=\sum_{n\geq 0}\frac{a^{n}}{n!}D^{n}$  is
the  translation operator  ($ D=\frac{d}{dx}$ ) and the recurrence
(\ref{C}) is obvious.\\

\noindent {\bf Example 2.} $ c_{n,k} =\Big\{ {n \atop k}\Big\}$  denote
Stirling numbers of the second type. Then
$x^{n}=\sum_{k=0}^{n}\Big\{ {n \atop k}\Big\}x^{\underline{k}}$ .
Let $B_{n}=\sum_{k=0}^{n}\Big\{ {n \atop k}\Big\}$ , then
\\$B_{n}= number\; of\; all\; partitions\; of \; S ;\; n =
|S|$;$B_{n}\equiv$ Bell numbers.

\noindent As the recurrence for Bell numbers reads $B_{n+1}=\sum_{k=0}^{n}\binom{n}{k}B_{k}\; ;\;\;n\geq 0$
   we easily get from (\ref{C}) the inspiring identity
$$\sum_{k=0}^{n}\binom{n}{k}B_{k}=B_{n}+\sum_{k=1}^{n}\Big\{ {n \atop k}\Big\}k.$$

\vspace{0.1cm}
\noindent {\bf Example 3.} $ c_{n,k} =\Big[
{n \atop k}\Big] $  are Stirling numbers of the first kind . Then
$x^{\overline{n}}=\sum_{k=0}^{n}\Big[ {n \atop k}\Big]
x^{k},\;n\geq 0$. $C_{n}=\sum_{k=0}^{n}\Big[ {n \atop k}\Big]=n!$
 and the obvious recurrence $C_{n+1}= ( n + 1) C_{n}$ coincides with (\ref{C}) .

\vspace{0.1cm}
\noindent {\bf Example 4.} $ c_{n,k} =\binom{n}{k}_{q} ,\; x^{n}=\sum_{k\geq
0}\binom{n}{k}_{q}\Phi_{k}(x),\;\;\;\Phi_{k}(x)=\prod_{s=0}^{k-1}(x-q^{s})$,
$$ C_{n}=\sum_{k \geq 0}\binom{n}{k}_{q} = number\;\; of\;\; all\;\;
subspaces\;\; of\;\; V(n,q)\;\; i.e.$$
$$C_{n}=\sum_{k=0}^{n}\binom{n}{k}_{q}=\left| L(n,q)\right|\equiv
(1+_{q}1)^{n},$$
where  $L(n,q)$ denotes the lattice of all subspaces of  $V(n,q)$ - ( the $n$-th dimensional space
$V(n,q)$  over Galois field   $GF(q)$  )  - see: \cite{4} .\\
According to Konvalina
$$C_{n}=\sum_{k=0}^{n}\binom{n}{k}_{q}\equiv (1+_{q}1)^{n}$$
is at the same time the number of all possible choices
\underline{with repetition} of  objects from exponentially
weighted boxes \cite{5,6}.\\
As $q$-Gaussian polynomials   $\Phi_{k}(x) =\prod_{s=0}^{k-1}(x-q^{s})$
constitute a persistent root polynomial sequence - the generalized
Lah numbers - which are now  are
$\binom{n}{k}_{q}=c_{n,k}=L_{n,k}$ must satisfy the recurrence
equation
$$L_{n+1,k}=L_{n,k-1}+r_{k+1}L_{n,k},\;\; L_{0,0}=1,\;\;L_{0,-1}=0\;\;\;n,k\geq 0,$$
where  $[r]=\left\{ q^{k-1}\right\}_{k\geq 1}$   is the root sequence determining
 $\left\{ \Phi_{n}\right\}_{n\geq 0}$. See Appendix 2.
Hence the recurrence   (\ref{C})  for the number $C_{n}$  of all
possible choices \underline{with  repetition} of  objects from
exponentially weighted boxes \cite{5,6}  is of the form
$$C_{n+1}=C_{n}+\sum_{k=0}^{n}\binom{n}{k}_{q}q^{k}$$
because
$$C_{n+1}=C_{n}+\sum_{k=0}^{n}c_{n,k}r_{k+1}.$$
In another   $q$-umbral notation (see: Appendix 2)
$$(1+_{q}1)^{n+1}=(1+_{q}1)^{n}+\sum_{k=0}^{n}\binom{n}{k}_{q}q^{k}.$$
Occasionally note that using the M\"{o}bius inversion formula
(see: \cite{4}  - the  $q$-binomial theorem Corollaries in section
5 - included )   we have
$$ x^{n}=\sum_{k\geq 0}\binom{n}{k}_{q}\Phi_{k}(x),\;\;\;\;
\Phi_{k}(x)=\sum_{l \geq
0}\binom{k}{l}_{q}(-1)^{l}q^{\binom{l}{2}}x^{k-l},\;\;n\geq 0$$
i.e. the identity
$$\prod_{s=0}^{k-1}(x-q^{s})\equiv \sum_{l\geq 0}\binom{k}{l}_{q}(_1)^{l}q^{\binom{l}{2}}x^{k-l}.$$
"The M\"{o}bius $q$-inverse example" is then the following.

\vspace{0.1cm}
\noindent {\bf Example 4'.}
$c_{n,k}=\binom{n}{k}_{q}(-1)^{n-k}q^{\binom{n-k}{2}} ,\;\;
x^{n}=\sum_{k\geq 0}c_{n,k}H_{k}(x)$, where (written with help of
generalized shift operator \cite{7,8,9,10} - analogously to the
use of $E^{-1}(D)x^{k}=(x-1)^{k}$ above in Example 1. ) we have
(see Appendix 2):
$$ H_{k}(x)=E(\partial_{q})x^{k}\equiv (x+_{q}1)^{k}=\sum_{l\geq
0}\binom{k}{l}_{q}x^{l}.$$
Compare with $H_{k}(x,y)|_{y=1}=H_{k}(x)$ in \cite{11}
(pp. 240-241) and note that
$$C_{n}=\sum_{k=0}^{n}|c_{n,k}|=\sum_{k=0}^{n}\binom{n}{k}_{q}q^{\binom{k}{2}}$$
is the number of all possible choices \underline{without
repetition} of objects from exponentially weighted boxes
\cite{5,6}. It is then combinatorial - natural to consider also
another triad .

\vspace{0.1cm}
\noindent {\bf Example 4''.}
$c_{n,k}=\binom{n}{k}_{q}q^{\binom{k}{2}},\;\; x^{n}=\sum_{k\geq
0}\binom{n}{k}_{q}q^{\binom{n}{2}}\Gamma_{k}(x),\\
\Gamma_{k}(x)=\sum_{l\geq
0}\binom{k}{l}_{q}(-1)^{k-l}x^{l},\;\;n\geq 0$.

\vspace{0.1cm}
\noindent {\bf Example 5.} This is The Fibonomial Example: $
c_{n,k}=\binom{n}{k}_{F},\\ x^{n}=\sum_{k\geq
0}\binom{n}{k}_{F}\Xi_{k}(x)$ where $\binom{n}{k}_{F}$
are called  Fibonomial coefficients (see A3).
Using the M\"{o}bius inversion formula - analogously to  the
Example 4 case - we get
$$\Xi_{k}(x)=\sum_{l\geq 0}\binom{k}{l}_{f}\mu(l,k)x^{k-1},\;\;n\geq 0,$$
 where  M\"{o}bius  matrix $\mu$  is the unique inverse
of the incidence matrix $\zeta$ defined for the "Fibonacci cobweb"
poset introduced in  \cite{12}  (see: A4 for explicit formula of
$\zeta$ from \cite{16} ). For combinatorial interpretation of
these Fibonomial coefficients see:  Appendix 4. Interpretation of
$C_{n}$ stated below results from combinatorial interpretation of
Fibonomial coefficients.  Combinatorial interpretation of  ccc
$C_{n}$ in this case is then the following:

$C_{n}=\sum_{k\geq 0}\binom{n}{k_{F}}=${\em number of
\underline{all} such subposets of the Fibonacci cobweb poset  $P$
which are cobweb subposets starting from the level  $F_{k}$ and
ending at the level labeled by $F_{n}$.}

\vspace{0.1cm}
\noindent As  in the $q$-umbral case of the Example 4  we have  "$F$-umbral"
representation of  the Fibonomial ccc  $C_{n}$   ( see A3) .
Namely
$$C_{n}=\sum_{k=0}^{n}\binom{n}{k}_{F}\equiv (1+_{F}1)^{n}$$

\section{Remark on ccc - recurrences "exercise"}

\vspace{0.1cm}
\noindent ccc-recurrences from   examples 1-4  were an easy game to play.
All these four cases of  ccc might be interpreted on equal footing
due to the  ingenious scheme of one unified combinatorial
interpretation by Konvalina \cite{5,6}.

\vspace{0.1cm}
\noindent  As for  the other cases   the "exercise"  of finding the
recurrences for ccc  - up to the Fibonomial  case -  might be the
slightly harder  task. In this last case  with  Fibonomial
coefficients     one may derive  the well known  recurrence
\begin{equation}\label{*_F}
 \binom{n+1}{k}_{F}=F_{k-1}\binom{n}{k}_{F}+F_{n-k+2}\binom{n}{k-1}_{F}
 \end{equation}
 $$\binom{n}{0}_{F}=1,\;\;\binom{0}{k}_{F}=0\;\; \mathrm{for}\;\;k>0$$
or equivalently
$$\binom{n+1}{k}_{F}=F_{k+1}\binom{n}{k}_{F}+F_{n-k}\binom{n}{k-1}_{F}$$
$$\binom{n}{0}_{F}=1,\;\;\binom{0}{k}_{F}=0\;\; \mathrm{for}\;\;k>0$$
also  by combinatorial reasoning  \cite{16}  (it might be a hard
exercise - not to check it but to prove it another way - a use of
\cite{2} - being recommended). Derivation of the recurrence for
Fibonomial ccc $C_{n}$ - we leave  as an exercise.

\section{Appendix}

\subsection{On the notation used trough out in this note}
$q$-binomial and  Fibonomial cases  are specifications of the
$\psi$-sequence choice in the so called  $\psi$- calculus notation
\cite{8}-\cite{10} in conformity with \cite{14} referring back  to
works of Brenke and Boas  on generalized Appell polynomials (see
abundant references in \cite{9}) \subsection{$q$-Gaussian ccc}
 \begin{equation}\label{*}
\binom{n+1}{k}_{q}=q^{k}\binom{n}{k}_{q}+\binom{n}{k-1}_{q}
 \end{equation}
 $$\binom{n}{0}_{q}=1,\;\;n\geq 0,k\geq 1,$$
 where $n_{q}!=n_{q}(n-1)_{q}!;\;\;1_{q}!=)_{q}!=1,\;\;n_{q}^{\underline{k}}=
 n_{q}(n-1)_{q}\ldots (n-k+1)_{q},\\ \binom{n}{k}_{q}\equiv
 \frac{n_{q}^{\underline{k}}}{k_{q}!}$.
 \begin{rem}{\em (see \cite{15}) The dual to (\ref{*}) reccurrenceis
 then given by
 \begin{equation}\label{**}
 x\Phi_{n}(x)=q^{n}\Phi_{n}(x)+\Phi_{n+1}(x)
 \end{equation}
 $$\Phi_{0}(x)=1,\;\;\Phi_{-1}(x)=0;\;\;n\geq 0,$$
 in accordance with a well known fact  (see: \cite{4} ) that
\begin{equation}\label{***}
x^{n}=\sum_{k\geq 0}\binom{n}{k}_{q}\Phi_{k}(x),
\end{equation}
where   $\Phi_{k}(x)=\prod_{s=0}^{k-1}(x-q^{s})$ are the  $q$-Gaussian polynomials.

\vspace{0.1cm}
\noindent Recall now that the Jackson $\partial_{q}$- derivative is defined as follows
\cite{7}-\cite{11}: $\left(
\partial_{q}\varphi\right)(x)=\frac{\varphi(x)-\varphi(qx)}{(1-q)x}$ - and is also called the
$q$-derivative. The consequent notation for the generalized shift
operator \cite{7}-\cite{10},\cite{13} is:
$$E^{a}(\partial_{q})x^{n}=\exp_{q}\{a\partial_{q}\}=\sum_{k=0}^{\infty}\frac{a^{k}}{k_{q}!}\partial_{q}^{k}.$$
 Then we identify $E^{y}(\partial_{q})x^{n}\equiv (x+_{q}y)^{n}$. (Ward
   in \cite{7} does not use any subscript like $q$  or $F$ or  $\psi$.
   
\vspace{0.2cm}
\noindent Naturally  his   $"+"$  means  not  $+$ but  a "plus with a
subscript" in our notation).}
 \end{rem}
 \subsection{Fibonomial  notation}
   In straightforward analogy consider now the {\em Fibonomial coefficients}
   $c_{n,k}$,
   $$c_{n,k}=\binom{n}{k}_{F}=\frac{F_{n}!}{F_{k}!F_{n-k}!}=\binom{n}{n-k}_{F},$$
    where   $n_{F}\equiv F_{n}$,$n_{F}!\equiv n_{F}(n-1)_{F}(n-2)_{F}(n-3)_{F}\ldots
2_{F}1_{F};\;\;0_{F}!=1$;$n_{F}^{\underline{k}}=n_{F}(n-1)_{F}\ldots
(n-k+1)_{F};\;\;\;\binom{n}{k}_{F}\equiv
\frac{n^{\underline{k}}_{F}}{k_{F}!}$ and linear difference
operator $\partial_{F}$  acting as follows:
$\partial_{F};\;\;\partial_{F}x^{n}=n_{F}x^{n-1};\;\;n\geq 0$ - we
shall call the $F$-derivative. Then in conformity with \cite{7}
and with notation as in \cite{10,9,8,13} one writes:
\renewcommand{\labelenumi}{(\arabic{enumi})}
\begin{enumerate}
\item $\left( x+_{F}a\right)^{n}\equiv \sum_{k\geq
0}\binom{n}{k}_{F}a^{k}x^{n-k}$  where  $\binom{n}{k}_{F}\equiv
\frac{n^{\underline{k}}_{F}}{k_{F}!}$ \\and
$n_{F}^{\underline{k}}=n_{F}(n-1)_{F}\ldots (n-k+1)_{F}$; \item
$\left( x+_{F}a\right)^{n}\equiv
E^{a}(\partial_{F}x^{n};\;\;E^{a}(\partial_{F})=\sum_{n\geq
0}\frac{a^{n}}{n_{F}!}\partial_{F}^{n};\;\;\\E^{a}(\partial_{F})f(x)=f(x+_{F}a),$
$E^{a}(\partial_{F})$ is corresponding generalized translation
operator.
\end{enumerate}

\subsection{Combinatorial interpretation of Fibonomial Coefficients
\cite{12,16,17,18,19,20} }

\vspace{0.1cm}
\noindent Let us depict a partially ordered infinite set  $P$ via
its finite part  Hasse diagram to be continued ad infinitum in an
obvious way as seen from the figure below.

\vspace{0.1cm}
\noindent This Fibonacci cobweb partially ordered infinite set  $P$ is defined as in \cite{12} via its
finite part - cobweb subposet  $P_{m}$({\em rooted at $F_{1}$
level subposet}) to be continued ad infinitum in an obvious way as
seen from the figure of  $P_{5}$  below.  It looks like the
Fibonacci tree with a specific "cobweb". It is identified with
Hasse diagram of the partial order set  $P$ .

\begin{figure}[ht]
\begin{center}
	\includegraphics[width=100mm]{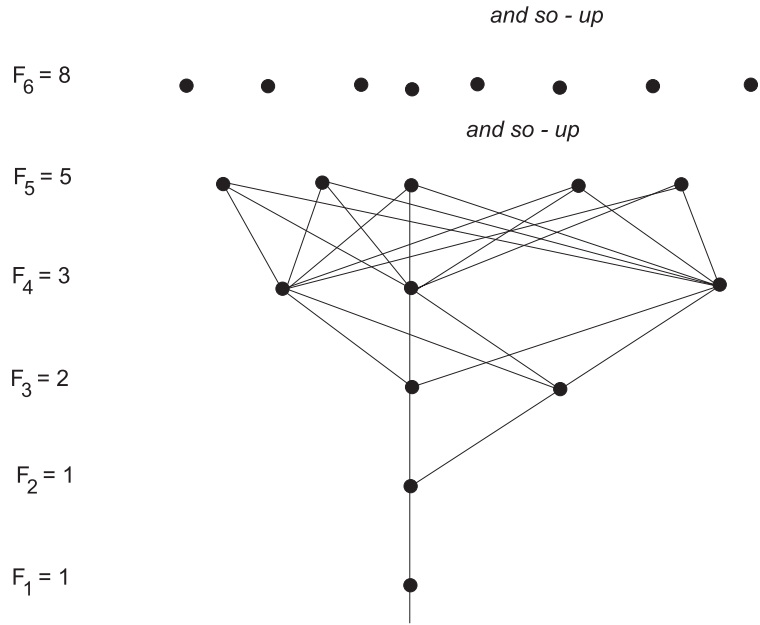}
	\caption { The graded structure of the  Fibonacci cobweb  poset}
\end{center}
\end{figure}

\vspace{0.1cm}
\noindent If one defines  this poset  $P$  with help of  its incidence
matrix $\zeta$ representing $P$ uniquely   then one arrives at
$\zeta$ with easily recognizable staircase-like  structure - of
zeros in the upper part of this upper triangle matrix   $\zeta$.
This structure is depicted by the Figure 2  where:  empty places
mean zero values (under diagonal)  and  in  filled with  - -...-   parts of rows
(above the diagonal) "`-"' stay for ones. The picture below is drawn for the sequence $F= \left\langle F_\textbf{ \textcolor{red}{1}}, F_2,F_3,...,F_n,... \right\rangle$, where 
$F_k$ are Fibonacci numbers.

\newpage

\begin{flushright}

1 - - - - - - - - - - - - - - - - - - - - - - - - - - - - - - - - - - - - - - - - - - - - - - - - - - -\\
1 - - - - - - - - - - - - - - - - - - - - - - - - - - - - - - - - - - - - - - - - - - - - - - - -\\
1\textbf{ \textcolor{red}{0}} - - - - - - - - - - - - - - - - - - - - - - - - - - - - - - - - - - - - - - - - - - - - - -\\
1 - - - - - - - - - - - - - - - - - - - - - - - - - - - - - - - - - - - - - - - - - - - - - - -\\
1 \textbf{\textcolor{red}{0 0}} - - - - - - - - - - - - - - - - - - - - - - - - - - - - - - - - - - - - - - - - - - -\\
1 \textbf{\textcolor{red}{0}} - - - - - - - - - - - - - - - - - - - - - - - - - - - - - - - - - - - - - - - - - - -\\
1 - - - - - - - - - - - - - - - - - - - - - - - - - - - - - - - - - - - - - - - - - - -\\$F_{5}-1\; \textbf{\textcolor{red}{0}}'s \;\; $1 \textbf{\textcolor{red}{0 0 0 0}} - - - - - - - - - - - - - - - - - - - - -  - - - - - - - - - - - - - - -\\
1 \textbf{\textcolor{red}{0 0 0}} - - - - - - - - - - - - - - - - - - - - - - - - - - - - - - - - - - - -\\
1 \textbf{\textcolor{red}{0 0}} - - - - - - - - - - - - - - - - - - - - - - - - - - - - - - - - - - - -\\
1 \textbf{\textcolor{red}{0}} - - - - - - - - - - - - - - - - - - - - - - - - - - - - - - - - - - - -\\
1 - - - - - - - - - - - - - - - - - - - - - - - - - - - - - - - - - - - -\\
$F_{6}-1\;zeros \quad \quad \quad \quad $1 \textbf{\textcolor{red}{0 0 0 0 0 0 0}} - - - - - - - - - - - - - - - - - - - - - - - - - -\\
1 \textbf{\textcolor{red}{0 0 0 0 0 0}} - - - - - - - - - - - - - - - - - - - - - - - - -\\
1 \textbf{\textcolor{red}{0 0 0 0 0}} - - - - - - - - - - - - - - - - - - - - - - - - -\\
1 \textbf{\textcolor{red}{0 0 0 0}} - - - - - - - - - - - - - - - - - - - - - - - - -\\
1 \textbf{\textcolor{red}{0 0 0}} - - - - - - - - - - - - - - - - - - - - - - - - -\\
1 \textbf{\textcolor{red}{0 0}} - - - - - - - - - - - - - - - - - - - - - - - - -\\
1 \textbf{\textcolor{red}{0}} - - - - - - - - - - - - - - - - - - - - - - - - -\\
1 - - - - - - - - - - - - - - - - - - - - - - - - -\\
$F_{7}-1\;zeros \quad \quad \quad \quad \quad \quad \quad \quad
\quad \quad   \;$1 \textbf{\textcolor{red}{\textcolor{red}{0 0 0 0 0 0 0 0 0 0 0 0}}} - - - - - - - \\
1\textbf{ \textcolor{red}{0 0 0 0 0 0 0 0 0 0 0}} - - - - - - -\\
1\textbf{ \textcolor{red}{0 0 0 0 0 0 0 0 0 0}} - - - - - - -\\
1 \textbf{\textcolor{red}{0 0 0 0 0 0 0 0 0}} - - - - - - -\\
1\textbf{ \textcolor{red}{0 0 0 0 0 0 0 0}} - - - - - - - \\
1\textbf{ \textcolor{red}{0 0 0 0 0 0 0}} - - - - - - - \\
1 \textbf{\textcolor{red}{0 0 0 0 0 0}} - - - - - - - \\
1 \textbf{\textcolor{red}{0 0 0 0 0}} - - - - - - - \\
1 \textbf{\textcolor{red}{0 0 0 0}} - - - - - - - \\
1 \textbf{\textcolor{red}{0 0 0}} - - - - - - - \\
1 \textbf{\textcolor{red}{0 0}} - - - - - - -  \\
1
\textbf{\textcolor{red}{0}} - - - - - - - \\
1 - - - - - - - \\
$F_{8}-1\;zeros \quad \quad \quad \quad \quad \quad \quad \quad
\quad \quad \quad \quad \quad \quad \quad \quad \quad \quad \quad \quad \;\,$.........................\\
{\em and so on}
\end{flushright}
\begin{center}
\noindent {\small Fig. \textbf{La Scala di Fibonacci}. The  staircase structure of incidence
matrix $\zeta \in I(P,R)$ = incidence algebra of $P$ over the commutative ring $R$. See [12](2003)} \end{center}

\newpage

\noindent For many pictures and further progress on cobweb posets and KoDAGs consult [17-29,31-42].

\vspace{0.2cm}

\noindent  The  characteristic function  $\chi_F(\leq)$ = $\zeta_F(\leq)$  of the partial order relation $\leq$ of the graded \textbf{$F$}-cobweb poset $P$ - or in brief just $\zeta$ - might be expressed by Kronecker delta  as here down  (\cite{18} (2003) and see also [12,40,18,16]) \textbf{for any natural numbers valued sequence $F$} - including the motivating example with  Fibonacci sequence $F$ denominated graded cobweb poset $P$ for which it reads ($F_0 = \textbf{\textcolor{red}{0}}$):   

\vspace{0.1cm}
\noindent  \textbf{Kwa\'sniewski $\zeta$ 2003 formula}
$$\zeta =\zeta_{1}-\zeta_{0}$$
where for $x,y \in {\bf N}$:
$$\zeta_{1}(x,y)=\sum_{k=0}^{\infty}\delta(x+k,y)$$
$$\zeta_{0}(x,y)=\sum_{k \geq 0}\sum_{s \geq 0}\delta
(x,F_{s+1}+k)\sum_{r=1}^{F_{s}-k-1}\delta (k+F_{s+1}+r,y)$$ and
naturally
$$ \delta (x,y)=\Big\{\begin{array}{l}1\;\;\;\;\; x=y\\0\;\;\;\;\;x\neq y\end{array}.$$

\vspace{0.2cm}

\noindent Right from the Hasse diagram of   $P$  here now obvious observations follow.

\vspace{0.2cm}

\begin{obs} The number of maximal chains starting from the root  (level  $F_{1}$ ) to reach any point at
the $n$-th level  labeled by  $F_{n}$  is equal to $n_{F}!$.
\end{obs}

\vspace{0.2cm}

\begin{obs}
The number of maximal chains starting from any fixed point at
the level labeled by  $F_{k}$  to reach any point at the $n$-th
level labeled by  $F_{n}$  is equal to
$n^{\underline{m}}_{F}\;\;(n = k+m )$ .
\end{obs}

\vspace{0.1cm}
\noindent Let us denote by $P_{m}$  a subposet of  $P$  with vertices
up to $m$-th level vertices 
$$\bigcup_{s=1}^{m}\Phi_{s}\;\; ;\;\Phi_{s} \;\;is\;\; the\;
\;set\;\;of\;\;elements\;\;of\;\;the\;\;s-th\;\;level$$.

\vspace{0.1cm}

\noindent \textbf{\textcolor{blue}{Attention please}.} $P_{m}$ depending on the context is also to be viewed on  as representing the set of maximal chains.
It is naturally a particular case of the layer of graded poset called cobweb poset [1-23] - for the reader convenience - see \textbf{\textcolor{red}{Appendix}} with basic ponderables. For the sake of coming right now combinatorial interpretation, all layers  on one hand are viewed on as  subposets while on the other hand are also to be viewed on  as representing the set of maximal chains.

\vspace{0.1cm}
\noindent  Consider now the following behavior of a {\em "sub-cob"}  moving from any
given point  of the  $F_{k}$  level of the poset  up. It behaves as it has been born right there and can reach at first $F_{2}$
points up then $F_{3}$ points up , $F_{4}$  and so on - thus climbing up to the level $F_{k+m} = F_{n}$ of the poset  $P$.  It
can see - as its Great Ancestor at the root $F_{1}$-th level  and potentially follow one of its own accessible finite subposet
$P_{m}$. One of many $P_{m}$ `s rooted at the $k$-th level might
be found. How many?

\vspace{0.2cm}

\noindent The answer \textbf{for the Fibonacci} sequence $F$ was given in  [12] (\textbf{2003}) :

\begin{obs}
 Let  $n = k+m $. The number of subposets equipotent to subposet
$P_{m}$ rooted at any \textbf{fixed}  point  at the level labeled by
$F_{k}$ and ending at the $n$-th level  labeled by  $F_{n}$ is
equal to
$$ \binom{n}{m}_{F}=\binom{n}{k}_{F}=\frac{n^{\underline{k}}_{F}}{k_{F}!}.$$
\end{obs}

\vspace{0.2cm}

\noindent Equivalently - now in a bit more mature \textbf{2009} year the  answer  is given simultaneously viewing layers as biunivoquely representing maximal chains sets. Let us make it formal.

\noindent Such recent equivalent formulation of this combinatorial interpretation is to be found in [20] from where we quote it here down. 

\noindent Let $\left\{ F_n \right\}_{n\geq 0}$ be a natural numbers valued sequence with $F_0 = 1$ (or  $F_0! \equiv 0!$ being exceptional as in case of Fibonacci numbers). Any such sequence uniquely designates both $F$-nomial coefficients of an $F$-extended umbral calculus as well as $F$-cobweb poset introduced  by this author (see :the source [19] from 2005 and earlier references therein). If these $F$-nomial coefficients are natural numbers or zero then we call the sequence $F$ - the $F$-\textbf{cobweb admissible sequence}.

\vspace{0.2cm}

\begin{defn}
Let any $F$-cobweb admissible sequence be given then $F$-nomial coefficients are defined as follows
$$
	\fnomial{n}{k} = \frac{n_F!}{k_F!(n-k)_F!} 
	= \frac{n_F\cdot(n-1)_F\cdot ...\cdot(n-k+1)_F}{1_F\cdot 2_F\cdot ... \cdot k_F}
	= \frac{n^{\underline{k}}_F}{k_F!}
$$
\noindent while $n,k\in N $ and $0_F! = n^{\underline{0}}_F = 1$.
\end{defn}

\vspace{0.2cm}

\begin{defn}

$C_{max}(P_n) \equiv  \left\{c=<x_0,x_1,...,x_n>, \: x_s \in \Phi_s, \:s=0,...,n \right\} $ i.e. $C_{max}(P_n)$ is the set of all maximal chains of $P_n$
\end{defn}

\vspace{0.2cm}

\begin{defn}
Let  $$C_{max}\langle\Phi_k \to \Phi_n \rangle \equiv \left\{c=<x_k,x_{k+1},...,x_n>, \: x_s \in \Phi_s, \:s=k,...,n \right\}.$$
Then the $C\langle\Phi_k \to \Phi_n \rangle $ set of  Hasse sub-diagram corresponding maximal chains defines biunivoquely 
the layer $\langle\Phi_k \to \Phi_n \rangle = \bigcup_{s=k}^n\Phi_s$  as the set of maximal chains' nodes and vice versa -
for  these \textbf{graded} DAGs (KoDAGs included).
\end{defn}

\vspace{0.2cm}

\noindent The \textbf{ equivalent} to that of of \textbf{Obsevation 3}  formulation of combinatorial interpretation of cobweb posets via their cover relation digraphs (Hasse diagrams) is the following.

\vspace{0.2cm}

\noindent \textbf{Theorem} [20]\\
\noindent(Kwa\'sniewski) \textit{For $F$-cobweb admissible sequences $F$-nomial coefficient $\fnomial{n}{k}$ is the cardinality of the family of \emph{equipotent} to  $C_{max}(P_m)$ mutually disjoint maximal chains sets, all together \textbf{partitioning } the set of maximal chains  $C_{max}\langle\Phi_{k+1} \to \Phi_n \rangle$  of the layer   $\langle\Phi_{k+1} \to \Phi_n \rangle$, where $m=n-k$.}

\vspace{0.2cm}

\vspace{0.2cm}
\noindent For February 2009 readings on further progress in combinatorial interpretation and application of the presented author invention  i.e. partial order posets named  cobweb posets and their's corresponding encoding Hasse diagrams KoDAGs see  [19-29,31-42] and references therein. For active presentation of cobweb posets see [40]. This electronic paper is an update of [41]

%%%%%%%%%%%%%%%%%%%%%%%%%%%%%%%%%%%%%%%%%%%%%%%%%%%%%%%%%%%%%%%%%%%%%%%%%%%%%%%%%%%%%%%%%%%%%%%%%%%%%%%%%%%%%%%%%%%%%%%%%55

\vspace{0.2cm}

\section{ Summarizing , Concluding  ad upgrade and ad references  remarks}

\vspace{0.2cm}

\noindent \textbf{Remark 1}
\noindent The matrix elements of $\zeta(x,y) $ were given in 2003 ([18] Kwa\'sniewski) using  $x,y\in N\cup\left\{0\right\} $  labels of  vertices in their "`natural"' linear order: \\
\textbf{1.} set  $k=0$,\\ 
\textbf{2.} then label subsequent vertices - from the left to the right - along the level  $k$,\\
\textbf{3.} repeat  2.  for $k \rightarrow  k+1$  until  $k=n+1$ ;  $n  \in N\cup\left\{\infty\right\}$ 

\vspace{0.2cm}
\noindent As the result we obtain  the $\zeta$ matrix for Fibonacci sequence as presented by the the Fig. \textit{La Scala di Fibonacci} dating back to 2003 [18,41].

\vspace{0.2cm}
\noindent  Inspired [8,9,13] by Gauss $n_q = q^0 + q^1 + ...+ q^{n-1}$ finite geometries numbers and in the spirit of Knuth "`notationlogy"' [30] we shall refer here also to the  upside down notation effectiveness as in [21-24,20]  or earlier in [8,9,10,13], (consult also the  Appendix in [33]). \textbf{As for} that upside down attitude  $F_n \equiv n_F$ being much more than "`just a convention"' to be used substantially in what follows as well as  for  the sake of completeness - \textbf{let us} quote it as The Principle according to  Kwas\'sniewski [42]
where it has been  formulated as an "`\textbf{of course}"'  Principle i.e. simultaneously  trivial and powerful statement.
\vspace{0.2cm}

\begin{quote}
\textbf{The Upside Down Notation Principle.} \\
\textbf{1. Let} the statement $s(F)$ depends only on the fact that $F$ is a natural numbers valued statement.\\
\textbf{2. Then} if one proves that $s(N)\equiv s(\left\langle n\right\rangle_{n \in N})$ is true  - the statement  $s(F)\equiv s(\left\langle n_F\right\rangle_{n \in N})$ is also true. Formally - use equivalence relation classes induced by co-images of $s : \left\{F\right\} \mapsto 2^{\left\{1\right\}}$ and proceed in a standard way.

\end{quote}

\vspace{0.2cm}

\vspace{0.2cm} 

\noindent In order to proceed further let us now recall-rewrite  purposely  here Kwa\'sniewski $2003$ - formula for $\zeta$ function of \textbf{arbitrary}  cobweb poset so as to see that its' algorithm rules authomatically  make it valid for all $F$-cobweb posets where $F$ is any natural numbers valued sequence  i.e. with \textcolor{red}{\textbf{$F_0 > 0$}}. $I(\Pi,R)$ stays for the incidence algebra of the poset $\Pi$ over the commutative ring $R$.

$$\zeta(x,y) =\zeta_{1}(x,y) -\zeta_{0}(x,y)$$

$$\zeta_{\textcolor{green}{\textbf{1}}}(x,y)=\sum_{k=\textcolor{blue}{\textbf{0}}}^{\infty}\delta(x+k,y)$$
$$\zeta_{0}(x,y)=\sum_{k \geq \textcolor{red}{\textbf{1}}}\sum_{s \geq 0}\delta
(x,F_{s+1}+k)\sum_{r=1}^{F_{s}-k-1}\delta (k+F_{s+1}+r,y)$$ and naturally
$$ \delta (x,y)=\Big\{\begin{array}{l}1\;\;\;\;\; x=y\\0\;\;\;\;\;x\neq y\end{array}.$$
\vspace{0.2cm}

\noindent The above formula  for  $\zeta \in I(\Pi,R)$ rewritten in ($F_s \equiv s_F$) upside down notation equivalent form as below
is of course \textbf{\textcolor{red}{valid for all}} cobweb posets.

$$\zeta(x,y) =\zeta_{\textcolor{green}{\textbf{\textcolor{green}{\textbf{1}}}}}(x,y) -\zeta_{0}(x,y),$$

$$\zeta_{\textcolor{green}{\textbf{1}}}(x,y)=\sum_{k=\textcolor{blue}{\textbf{0}}}^{\infty}\delta(x+k,y),$$ 
$$\zeta_{0}(x,y)=\sum_{s \geq 1}\sum_{k \geq \textcolor{red}{\textbf{1}}}\delta (x,k+s_F)\sum_{r=1}^{(s-1)_F-k-1}\delta (x+r,y).$$ \\

\textbf{Note.} $+\zeta_{\textcolor{green}{\textbf{1}}}$ "`produces the Pacific ocean of  \textcolor{green}{1's}"' in the whole upper triangle part of a would be incidence algebra $\sigma \in I(\Pi,R)$ matrix elements\\

\textbf{Note.}  $-\zeta_0$ cuts out \textcolor{red}{0's} i.e.  thus producing "`zeros' $F$-La Scala staircase"' in the \textcolor{green}{\textbf{1}'s}  delivered by $+\zeta_{\textcolor{green}{\textbf{1}}}$.\\

\vspace{0.2cm}

\noindent  This results exactly in forming \textcolor{red}{0's} rectangular triangles:  $s_F - 1$ of them at the start of subsequent stair and then  down to one \textcolor{red}{0} till - after  $s_F - 1$ 
rows passed by one reaches a half-line of $1's$ which is  running to the right- right  to infinity and thus marks the next in order stair of the $F$- La Scala.\\

\vspace{0.2cm}

\noindent The $\zeta$ matrix explicit formula was given for arbitrary graded posets with the finite set of minimal in terms of natural join of bipartite digraphs in SNACK = the \textbf{S}ylvester \textbf{N}ight \textbf{A}rticle on\textbf{ K}oDAGs and\textbf{ C}obwebs = [21].\\

\vspace{0.2cm}

\noindent What was said is equivalent to the fact that the cobweb poset coding La Scala is of the natural join operation origin while thus producing  $\zeta$ matrix [23,21,24] which is of the form 
(quote from SNACK = [21], see: Subsection 2.6.)

\vspace{0.2cm}

\begin{quote}
\noindent The explicit  expression for zeta matrix $\zeta_F$ of cobweb posets  via known blocks of zeros and ones for arbitrary natural numbers valued $F$- sequence  was given in (here) [23]  due to more than  mnemonic  efficiency  of the up-side-down notation being applied (see [23] and references therein). With this notation inspired by Gauss  and replacing  $k$ - natural numbers with   "$k_F$"  numbers one gets 
$$
	\mathbf{A}_F = \left[\begin{array}{llllll}
	0_{1_F\times 1_F} & I(1_F \times 2_F) & 0_{1_F \times \infty} \\
	0_{2_F\times 1_F} & 0_{2_F\times 2_F} & I(2_F \times 3_F) & 0_{2_F \times \infty} \\
	0_{3_F\times 1_F} & 0_{3_F\times 2_F} & 0_{3_F\times 3_F} & I(3_F \times 4_F) & 0_{3_F \times \infty} \\
	0_{4_F\times 1_F} & 0_{4_F\times 2_F} & 0_{4_F\times 3_F} & 0_{4_F\times 4_F} & I(4_F \times 5_F) & 0_{4_F \times \infty} \\
	... & etc & ... & and\ so\ on & ...
	\end{array}\right]
$$

\noindent and

$$
	\zeta_F = exp_\copyright[\mathbf{A}_F] \equiv (1 - \mathbf{A}_F)^{-1\copyright} \equiv I_{\infty\times\infty} + \mathbf{A}_F + \mathbf{A}_F^{\copyright 2} + ... =
$$
$$
	= \left[\begin{array}{lllll}
	\textcolor{red}{\textbf{I}}_{1_F\times 1_F} & I(1_F\times\infty) \\
	O_{2_F\times 1_F} & \textcolor{red}{\textbf{I}}_{2_F\times 2_F} & I(2_F\times\infty) \\
	O_{3_F\times 1_F} & O_{3_F\times 2_F} & \textcolor{red}{\textbf{I}}_{3_F\times 3_F} & I(3_F\times\infty) \\
	O_{4_F\times 1_F} & O_{4_F\times 2_F} & O_{4_F\times 3_F} & \textcolor{red}{\textbf{I}}_{4_F\times 4_F} & I(4_F\times\infty) \\
	... & etc & ... & and\ so\ on & ...
	\end{array}\right]
$$

\noindent where  $I (s\times k)$  stays for $(s\times k)$  matrix  of  ones  i.e.  $[ I (s\times k) ]_{ij} = 1$;  $1 \leq i \leq  s,  1\leq j  \leq k.$  and  $n \in N \cup \{\infty\}$
\end{quote}

\vspace{0.1cm}

\noindent In the $\zeta_F $ formula from [23]  $\copyright$ denotes the Boolean  product, hence - used for Boolean powers too. We readily recognize from its block structure that $F$-La Scala 
is formed by \textcolor{red}{\textbf{upper zeros}} of block-diagonal matrices $\textcolor{red}{\textbf{I}}_{k_F\times k_F}$ which sacrifice   these their  \textcolor{red}{\textbf{zeros}} to constitute the  $k$-th   subsequent stair in the $F$-La Scala descending and descending far away down
to infinity. Thus the cobweb poset coding La Scala is due to the natural join origin of $\zeta$ matrix.

\vspace{0.2cm}

\noindent  Note now that because of $\delta$'s  under summations in the former $\zeta$ formula  the following is obvious:  

$$ 1 \leq r = y - x \leq (s-1)_F - k - 1 \equiv  1 \leq r = y - k - s_F \leq s-1)_F - k -1  \equiv $$
$$ \equiv  1 \leq r = y \leq s_F - (s-1)_F  - 1 .  $$

\vspace{0.1cm}

\noindent Because of that the above last expression of the $\zeta$ expressed in terms of  $\delta \in I(\Pi,R)$  may be still simplified [for the sake of verification and portraying via computer simple program implementation]. Namely the following is true:

$$\zeta(x,y) =\zeta_{\textcolor{green}{\textbf{1}}}(x,y) -\zeta_{0}(x,y),$$

\vspace{0.1cm}

\noindent where

$$\zeta_{\textcolor{green}{1}}(x,y)=\sum_{\textcolor{blue}{k=\textbf{0}}}^{\infty}\delta(x+k,y),$$ 

[- note: $+\zeta_1$ "`produces the Pacific ocean of \textcolor{green}{1's}"' in the whole upper triangle part of a would be incidence algebra $\sigma \in I(\Pi,R)$ matrix elements],\\

\vspace{0.1cm}

\noindent and where

$$\zeta_{0}(x,y)=\sum_{s \geq 1}\sum_{k \geq \textcolor{red}{\textbf{1}}}\delta (x,k+s_F)\sum_{r\geq 1}^{s_F + (s-1)_F-1} \delta (r,y),$$ 

\noindent [- note then again that $-\zeta_0$ cuts out "one's $F$-La Scala staircase"' in the \textcolor{green}{1's} provided by $+\zeta_{\textcolor{green}{\textbf{1}}}$].\\

\vspace{0.1cm}

\noindent Note, that for \textcolor{red}{\textbf{$F$ = Fibonacci}} this still more simplifies   as then  $$s_F + (s-1)_F-1 = (s+1)_F.$$

\vspace{0.1cm}

\vspace{0.2cm}

\noindent \textbf{Remark 1.1. ad Knuth notation [30] indicated back to me by Maciej Dziemia\'nczuk}\\

\vspace{0.1cm} 

\noindent In the  wise "`notationlogy"' note [30]  among others one finds notation just for the purpose here (Maciej Dziemia\'nczuk's observation)

$$ [ s ] =\left\{ \begin{array}{cl} 1&if \ s \ is \  true ,\\0&otherwise.

\end{array} \right. $$

\noindent Using this makes my last above expression of the $\zeta$ in terms of  $\delta$ still more transparent and handy  if rewritten in Donald Ervin Knuth's  notation [30]. Namely:

$$\zeta(x,y) =\zeta_{\textcolor{green}{\textbf{1}}}(x,y) -\zeta_{0}(x,y)$$

$$\zeta_{\textcolor{green}{1}}(x,y)= [x \leq y]$$ 

$$\zeta_{0}(x,y)=\sum_{s \geq 1}\sum_{k \geq \textcolor{red}{\textbf{1}}}[x=k+s_F][1 \leq y \leq s_F + (s-1)_F\! - 1 ].$$

\vspace{0.1cm}

\noindent Note, that for \textcolor{red}{\textbf{$F$ = Fibonacci}} this still more simplifies   as then  $$s_F + (s-1)_F-1 = (s+1)_F.$$

\vspace{0.2cm}

\noindent \textbf{Remark 1.2. Knuth notation [30] - and  Dziemia\'nczuk's guess}

\vspace{0.2cm}

\noindent It was remarked by my \textbf{Gda\'nsk University} Student \textbf{Maciej Dziemia\'nczuk} - that my  $\zeta \in I(\Pi,R)$  (equivalent) expressions are valid according to him only for $F$ = Fibonacci sequence.  In view of the Upside Down Notation Principle if any of these  is  valid for any particular natural numbers  valued sequence $F$ it should be true for all of the kind.

\vspace{0.2cm}

\noindent His being doubtful  led him to creative invention of his own. 

\vspace{0.1cm}

\noindent  Here comes the  formula  postulated by him.

\vspace{0.1cm}

$$ \zeta(x,y) =  [x\leq y] - [x<y] \sum_{n\geq 0} [(x> S(n)] [y \leq S(n+1)],$$

\noindent where 

$$ S(n) = \sum_{k\geq 1}^n k_F $$

\vspace{0.1cm}

\noindent \textbf{Exercise.} My reply to this guess is the following  Exercise.

\noindent  Let $x,y\in N\cup\left\{0\right\}$  be the labels of  vertices in their "`natural"' linear order as explained earlier. 

\noindent Prove the claim:\\
\textit{Dziemia\'nczuk formula is equivalent to Kwa\'sniewski formulas}.\\ 

\noindent - What is against? \\
- What is for?  My "`for"' is the Socratic Method  question.  Why not use  the arguments in favor of 
$$\zeta_{0}(x,y)=\sum_{s \geq 1}\sum_{k \geq \textcolor{red}{\textbf{1}}}\delta (x,k+s_F)\sum_{r\geq 1}^{s_F + (s-1)_F-1} \delta (r,y),$$ 

\vspace{0.1cm}

\noindent  Hint. Use the same argumentation.
\vspace{0.1cm}

\noindent Here are some illustrative examples-exercises with pictures [Figures 2,3,4]  delivered by Maciej Dziemia\'nczuk's computer personal service using the above  Dziemia\'nczuk formula.

\begin{figure}[ht]
\begin{center}
	\includegraphics[width=100mm]{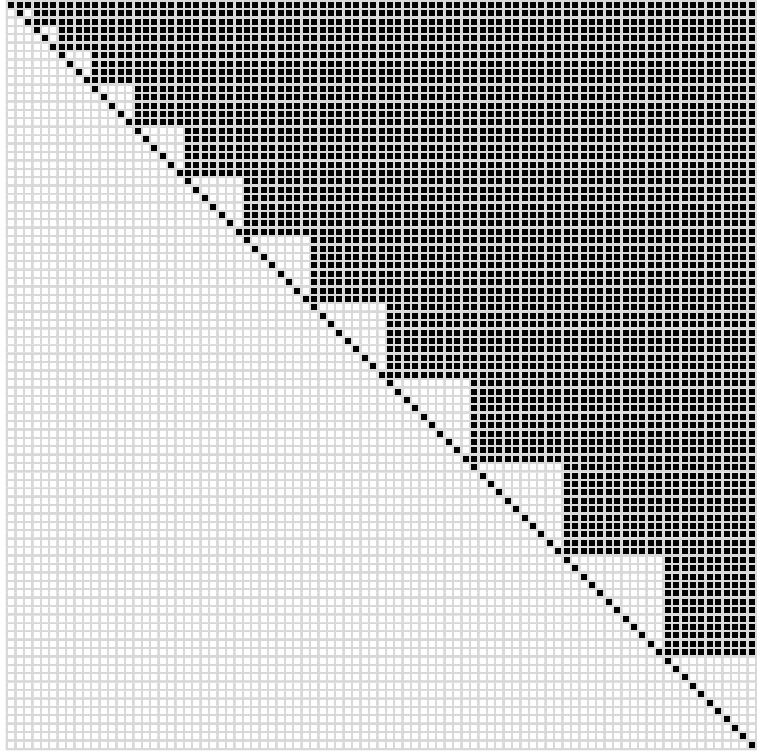}
	\caption {Display of the  $\zeta$ = $90 \times 90)$. The subposet $P_t$ of the  $N$ i.e. integers sequence $N$-cobweb poset. $t = ?$}
\end{center}
\end{figure}

\begin{figure}[ht]
\begin{center}
	\includegraphics[width=100mm]{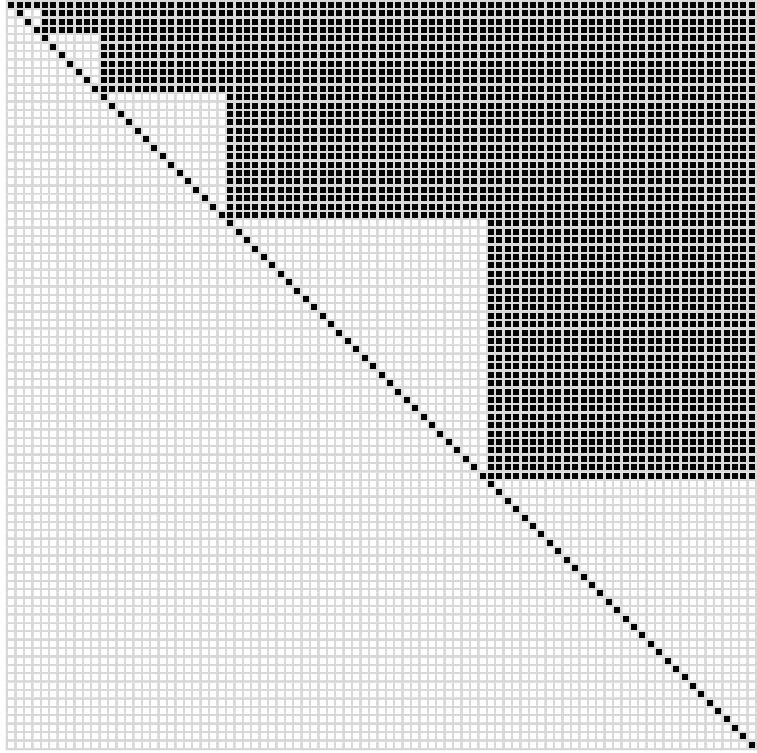}
	\caption {Display of the  $\zeta$ = $90 \times 90)$. The subposet $P_t$ of the  $F$ = Gaussian integers sequence $(q=2)$.  $F$-cobweb poset.$t = ?$}
\end{center}
\end{figure}

\begin{figure}[ht]
\begin{center}
	\includegraphics[width=100mm]{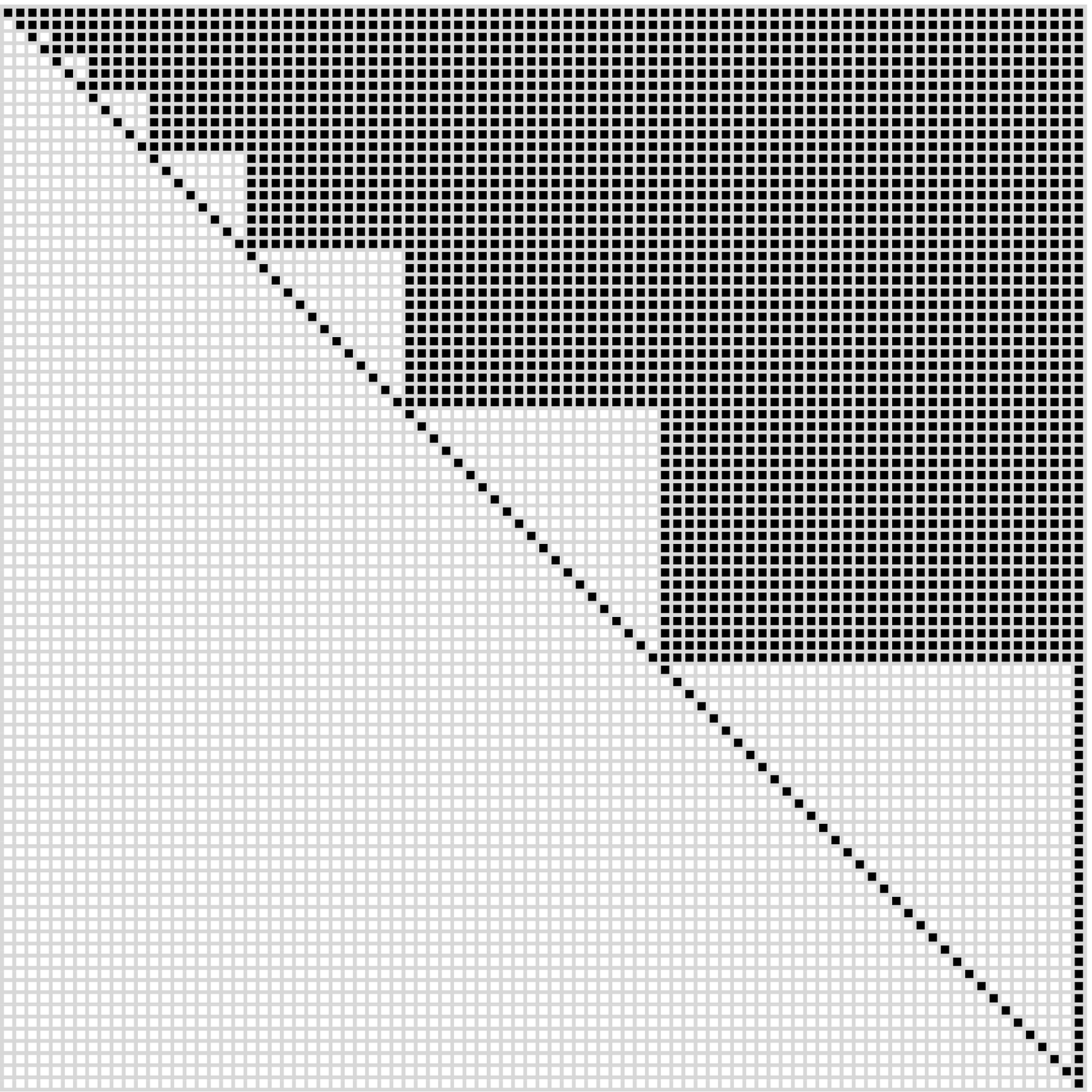}
	\caption {Display of the  $\zeta$ = $90 \times 90)$. The subposet $P_t$ of the  $F$ = Fibonacci-cobweb poset.$t = ?$}
\end{center}
\end{figure}

\vspace{0.2cm}

\noindent \textbf{Remark 2. Krot Choice.} While the above is established it is a matter of  simple observation by inspection  to find out how does  the the M{\"{o}}bius matrix $\mu =\zeta^{-1} $ looks like. Using in [25,26] this author example and expression for $\zeta$ matrix this  has been accomplished first (see also [26])   for Fibonacci sequence and then - via automatic extension - the same formula was given for $F$ sequences as above  by my  former student in her recent articles [27,28].  Here is her formula for  the cobweb posets' M{\"{o}}bius function (see: (6) in [28]).\\

$$\mu(x,y) = \mu (\left\langle s,t\right\rangle,\left\langle u,v\right\rangle) = \delta(s,u)\delta(t,v) - \delta(t+1,v) +  \sum_{k=2}^{\infty}\delta(t+k,v)(-1)^k\prod_{i= t+1}^{v-1}F_i-1 $$

\vspace{0.2cm}

\noindent Now bearing in mind the Upside Down Notation Principle  \textbf{for all} $F$-cobweb posets with   $F_0>0$ (as it should be for natural numbers valued sequences)  we may now rewrite the above in coordinate grid  $Z \times Z$ as below.

\vspace{0.2cm}
\noindent Let $x= \left\langle s,t\right\rangle$ and  $y= \left\langle u,v\right\rangle$ where $1 \leq s \leq F_t$, $1 \leq u \leq F_v$ while $t,v \in N$. Then

$$\mu(x,y) = \mu (\left\langle s,t\right\rangle,\left\langle u,v\right\rangle) = \delta(s,u)\delta(t,v)  +  \sum_{k=1}^{\infty}\delta(t+k,v)(-1)^k\prod_{i= t+1}^{v-1}(i-1)_F $$

\vspace{0.2cm}

\noindent The further relevant references of the same author see [27-29]. The above rewritten M{\"{o}}bius function formula is of course literally \textbf{valid for all} natural numbers valued sequences $F$. Consult: the  recent note "`On Characteristic Polynomials of the Family of Cobweb Posets"' [29].

\vspace{0.2cm}

\noindent The author of  [25] introduces parallely also another form of  $\zeta$ function formula and since now on -except for [27]- in  subsequent papers [26,28,29] their author  uses the formula for  $\zeta$  function in this another form. Namely - this  other form formula for $\zeta$  function in the present authors' grid coordinate system description of the cobweb posets was given by  Krot in her  note on M{\"{o}}bius function and M{\"{o}}bius inversion formula for Fibonacci cobweb poset [23]  with $F$ designating the Fibonacci cobweb  posets. In [24] the  formula for the  M{\"{o}}bius function for Fibonacci sequence $F$ was rightly treated as literally \textbf{valid for all }natural numbers valued sequences $F$ and the alternative. Of course the same concerns the $\zeta$  function formulas - the former and the latter. Here is this other latter form   formula for  $\zeta$  function (see: (7) in [25] or (1) in [27]).\\

\noindent Let $x= \left\langle s,t\right\rangle$ and  $y= \left\langle u,v\right\rangle$ where $1 \leq s \leq F_t$, $1 \leq u \leq F_v$ while $t,v \in N$. Then

$$\zeta (x,y) = \zeta (\left\langle s,t\right\rangle,\left\langle u,v\right\rangle) = \delta(s,u)\delta(t,v) +  \sum_{k=1}^{\infty}\delta(t+k,v) $$

\noindent where here - recall  ($a,b \in Z$): 

$$   \delta(a,b)=\left\{ \begin{array}{cl} 1&for\ a = b,\\0&otherwise.

\end{array} \right. $$

\vspace{0.2cm}

\noindent \textbf{Farewell interactive question.}  \textit{Are then all the presented $\zeta$ incidence function matrix of $F$- denominated cobweb posets equivalent?}

\section {Appendix [42]}

\vspace{0.2cm}

\noindent \textbf{Cobweb posets and  KoDAGs' ponderables  of the authors relevant productions.}

\begin{defn}
\noindent  Let  $n\in N \cup \left\{0\right\}\cup \left\{\infty\right\}$. Let   $r,s \in N \cup \left\{0\right\}$.  Let  $\Pi_n$ be the graded partial ordered set (poset) i.e. $\Pi_n = (\Phi_n,\leq)= ( \bigcup_{k=0}^n \Phi_k ,\leq)$ and $\left\langle \Phi_k \right\rangle_{k=0}^n$ constitutes ordered partition of $\Pi_n$. A graded poset   $\Pi_n$  with finite set of minimal 
elements is called \textbf{cobweb poset} \textsl{iff}  
$$\forall x,y \in \Phi \  i.e. \  x \in \Phi_r \ and \  y \in \Phi_s \   r \neq s\ \Rightarrow \   x\leq y   \ or \ y\leq x  , $$ 
 $\Pi_\infty \equiv \Pi. $
\end{defn}

\vspace{0.1cm}

\noindent \textbf{Note}. By definition of $\Pi$ being graded its  levels    $\Phi_r \in \left\{\Phi_k\right\}_k^\infty$ are independence sets  and of course partial order  $\leq $ up there in Definition 6.1. might be replaced by $<$.

\vspace{0.2cm}

\noindent The Definition 6.1.  is the reason for calling Hasse digraph $D = \left\langle \Phi, \leq \cdot \right\rangle $ of the poset $(\Phi,\leq))$ a \textbf{\textcolor{red}{Ko}}DAG as in  Professor   
\textbf{\textcolor{red}{K}}azimierz   \textbf{\textcolor{red}{K}}uratowski native language one word \textbf{\textcolor{red}{Ko}mplet} means \textbf{complete ensemble} - see more in  [23]
and for the history of this name see:  The Internet Gian-Carlo Polish Seminar \textbf{Subject 1.  oDAGs and KoDAGs in Company} (Dec. 2008).

\begin{defn}
\noindent Let  $F = \left\langle k_F \right\rangle_{k=0}^n$ be an arbitrary natural numbers valued sequence, where $n\in N \cup \left\{0\right\}\cup \left\{\infty\right\}$. We say that the cobweb poset $\Pi = (\Phi,\leq)$ is \textcolor{red}{\textbf{denominated}} (encoded=labelled) by  $F$  iff   $\left|\Phi_k\right| = k_F$ for $k = 0,1,..., n.$
\end{defn}

% % % % % % % % % % % % % % % % % % % % % % % % % % % % % % % % % % % % % % % % % % % % % % % % % % % % % % % % % % % % % % % % % % % % 

\vspace{0.4cm}

\noindent \textbf{Acknowledgments}
\vspace{0.1cm}
\noindent Thanks are expressed here to the now  Student of Gda\'nsk University Maciej Dziemia\'nczuk for applying his skillful   TeX-nology with respect to most of my articles since three years as well as for his general assistance and cooperation on KoDAGs  investigation.  Maciej Dziemia\'nczuk was not allowed to write his diploma with me being supervisor - while Maciej  studied in the local Bialystok University where my professorship till 2009-09-30  comes from.

\end{document}